\documentclass[a4paper,12pt]{article}
\usepackage{comment}
\usepackage{cite}
\usepackage{amsmath}
\usepackage{amssymb}
\usepackage{amsfonts}
\usepackage[T1]{fontenc}
\usepackage[utf8]{inputenc}
\usepackage{graphicx}
\usepackage{fancyhdr}
\usepackage{float}
\usepackage{xcolor}
\usepackage{authblk}
\usepackage{mathrsfs}
\usepackage{empheq}
\usepackage{url}

\usepackage{geometry}
\DeclareMathOperator{\sech}{sech}
\DeclareMathOperator{\csch}{csch}
\DeclareMathOperator{\arccosh}{arcCosh}
\DeclareMathOperator{\arcsinh}{arcSinh}
\DeclareMathOperator{\arctanh}{arcTanh}

\begin{document}

\title{A double series for $\pi$ using Fourier series and the Grothendieck-Krivine constant}

\author[$\dagger$]{Jean-Christophe {\sc Pain}\footnote{jean-christophe.pain@cea.fr}\\
\small
CEA, DAM, DIF, F-91297 Arpajon, France\\
Universit\'e Paris-Saclay, CEA, Laboratoire Mati\`ere en Conditions Extr\^emes,\\ 
91680 Bruy\`eres-le-Ch\^atel, France
}

\maketitle

\begin{abstract}
We provide a double-series formula for $\pi$ obtained using the Fourier series expansion of $1/\cos\left(x/4\right)$ and applying the Parseval-Plancherel identity. We show that such a formula involves the Grotehendieck-Krivine constant, and that the latter can therefore be expressed as a double series as well. 
\end{abstract}

\section{Introduction}

The seek for series expressions of $\pi$ has challenged mathematicians for many centuries (see for instance the Leibniz-Gregory-Madhava, Nilakantha, Ramanujan and Bailey-Borwein-Plouffe series) \cite{Borwein}. The purpose can be either to obtain the best compromise between simplicity and accuracy, to extract one specific digit without knowing the others, or to exhibit connections between $\pi$ and other mathematical quantities such as the golden ratio, the constants named after Catalan, Euler, Theodorus, {\it etc.}

An important contribution of Grothendieck to Banach space theory is the fundamental theorem in the metric theory of tensor products \cite{Grothendieck1956}. This theorem, also referred to as ``Grothendieck's inequality'', had a major impact in $C^{\ast}$-algebra theory \cite{Pisier1978}. A version of the theorem has been developed in the framework of ``operator spaces'' or non-commutative Banach spaces. Grothendieck's inequality is strongly related, in quantum mechanics, to the Bell inequalities \cite{Bell1964} for probabilities in classical locally deterministic theory and probabilities that arise in quantum theory. Grothendieck's inequality was crucial to test the Einstein-Podolsky-Rosen (EPR) framework of ``hidden variables'' proposed as an alternative to quantum mechanics \cite{Einstein1935}. Using Bell's ideas, experiments were carried out (see Refs. \cite{Aspect2000a,Aspect2000b}) to verify the presence of a deviation that invalidated the EPR conception. Tsirelson observed that the Grothendieck constant could be interpreted as an upper bound for the deviation in the generalized Bell inequalities (in particular there would be no deviation if the Grothendieck constant was equal to 1) \cite{Tsirelson1980,Tsirelson1993}. The Grothendieck constant has been introduced in graph theory and in computer science where the Grothendieck inequality is invoked to replace certain NP problems by others that can be treated by ``semidefinite programming'' and hence solved in polynomial time \cite{Pisier2012}.

In 1977, Krivine introduced the constant
\begin{equation}\label{KG}
K_G=\frac{\pi}{2\ln\left(1+\sqrt{2}\right)}, 	
\end{equation}
which he called ``Grothendieck's constant'' \cite{Krivine1977,Lelionnais1983}, but a better upper bound for Grothendieck's inequality was published in 2011 by Braverman, Makarychev, Makarychev and Naor  \cite{Braverman2011,Braverman2013}. Therefore, in order to avoid confusion, the constant given by formula (\ref{KG}) will be named the Grothendieck-Krivine in the following.

In section  \ref{sec2}, we present the Fourier series of a specific function and show that the Parseval-Plancherel enables one to obtain a double series formula for $\pi$. In section \ref{sec3}, we mention that the latter formulas involves the Grothendieck-Krivine constant, and provides therefore a connection between the latter and $\pi$ {\it via} a simple double series.

\section{Fourier series of $x\longmapsto \cos\left(\frac{x}{4}\right)$}\label{sec2}

Let us consider the function 
\begin{equation}
\mathscr{F}(x)=\frac{1}{\cos\left(\frac{x}{4}\right)}
\end{equation}
defined on $[-\pi,\pi]$ Its Fourier series reads
\begin{equation}
\mathscr{F}(x)=\frac{a_0}{2}+\sum_{n=1}^{\infty}a_n\cos(nx)+\sum_{n=1}^{\infty}b_n\sin(nx).
\end{equation}
Since the function is even, all the $b_n$ coefficients are equal to zero. We have in particular \cite{Tolstov1962,Kaplan1992}:
\begin{equation}
a_n=\frac{2}{\pi}\int_0^{\pi}\mathscr{F}(t)\cos(nt)dt=\frac{8}{\pi}\int_0^{\pi/4}\frac{\cos(4nt)}{\cos(t)}dt
\end{equation}
and
\begin{equation}
a_0=\frac{8}{\pi}\int_0^{\pi/4}\frac{dt}{\cos t}=\frac{8}{\pi}\ln\left[\tan\left(\frac{3\pi}{8}\right)\right].
\end{equation}
Taking into account the fact that
\begin{equation}
\tan\left(\frac{3\pi}{8}\right)=\frac{1-\cos\left(\frac{3\pi}{4}\right)}{\sin\left(\frac{3\pi}{4}\right)}=1+\sqrt{2},
\end{equation}
we get
\begin{equation}
a_0=\frac{8}{\pi}\ln(\sqrt{2}+1).
\end{equation}
Let us now calculate the difference of two successive coefficients $a_n-a_{n-1}$ \cite{Tissier1991}:
\begin{eqnarray}
a_n-a_{n-1}&=&\frac{8}{\pi}\int_0^{\pi/4}\frac{-2\left[\sin(4n-2)t\right]\sin(2t)}{\cos t}dt=-\frac{32}{\pi}\int_0^{\pi/4}\sin\left[(4n-2)t\right]\sin tdt\nonumber\\
& &-\frac{16}{\pi}\int_0^{\pi/4}\left\{\cos\left[(4n-3)t\right]-\cos\left[(4n-1)t\right]\right\}dt\nonumber\\
&=&-\frac{16}{\pi}\left\{\frac{\sin\left[(4n-3)\frac{\pi}{4}\right]}{4n-3}-\frac{\sin\left[(4n-1)\frac{\pi}{4}\right]}{4n-1}\right\}\nonumber\\
&=&\frac{16}{\pi}(-1)^n\left[\frac{\sin\left(\frac{\pi}{4}\right)}{4n-3}-\frac{\sin\left(\frac{3\pi}{4}\right)}{4n-1}\right]\nonumber\\
&=&\frac{8\sqrt{2}}{\pi}(-1)^n\left[\frac{1}{4n-3}-\frac{1}{4n-1}\right].
\end{eqnarray}
Therefore, one has
\begin{equation}
a_n=\frac{8}{\pi}\ln(1+\sqrt{2})-\frac{8\sqrt{2}}{\pi}\mathscr{S}_n
\end{equation}
where
\begin{equation}
\mathscr{S}_n=\sum_{k=1}^n(-1)^{k-1}\left(\frac{1}{4k-3}-\frac{1}{4k-1}\right).
\end{equation}
Since $\lim_{n\rightarrow\infty}a_n=0$, 

\begin{equation}
\mathscr{S}=\lim_{n\rightarrow\infty}\mathscr{S}_n=\frac{\sqrt{2}}{2}\ln(1+\sqrt{2})    
\end{equation}
and we have also
\begin{equation}
a_n=\frac{8\sqrt{2}}{\pi}\left(\mathscr{S}-\mathscr{S}_n\right).
\end{equation}
Thus
\begin{equation}
a_n=\frac{8\sqrt{2}}{\pi}\sum_{k=n+1}^{\infty}(-1)^{k}\left(\frac{1}{4k-1}-\frac{1}{4k-3}\right).
\end{equation}
The Parseval-Plancherel theorem \cite{Parseval1806,Plancherel1910} states that
\begin{equation}
\int_{-\pi}^{\pi}\frac{dt}{\cos^2\left(\frac{t}{4}\right)}=8\int_{0}^{\pi/4}\frac{du}{\cos^2\left(u\right)}=8
\end{equation}
and finally
\begin{empheq}[box=\fbox]{align}
\sum_{n=1}^{\infty}\left[\sum_{k=n+1}^{\infty}(-1)^{k}\left(\frac{1}{4k-1}-\frac{1}{4k-3}\right)\right]^2=\frac{\pi}{16}-\frac{1}{4}\ln^2\left(\sqrt{2}+1\right).
\end{empheq}

\section{The Grothendieck-Krivine constant}\label{sec3}

Let ${\bf A}$ be an $n\times n$ real square matrix with $n\ge 2$ such that \cite{mathe}:
\begin{equation}
\left|\sum_{i=1}^n\sum_{j=1}^na_{ij}s_it_j\right|\leq 1 	
\end{equation}
for all real numbers $s_1$, $s_2$, ..., $s_n$ and $t_1$, $t_2$, ..., $t_n$ satisfying $\left|s_i\right|$, $\left|t_j\right|\leq 1$. Then Grothendieck showed that there exists a constant $K_R(n)$ ensuring
\begin{equation}
\left|\sum_{i=1}^n\sum_{j=1}^na_{ij}x_i.y_j\right|\leq K_R(n) 	
\end{equation}
for all vectors $x_1$, $x_2$,..., $x_m$ and $y_1$, $y_2$,..., $y_n$ in a Hilbert space with norms $\left|x_i\right|\leq 1$ and $\left|y_j\right|\leq 1$. The Grothendieck constant is the smallest possible value of $K_R(n)$. As mentioned in the introduction, Krivine postulated that the limit 
\begin{equation}
\lim_{n\rightarrow \infty}K_R(n).	
\end{equation}
is equal to \cite{Krivine1977}:
\begin{equation}\label{KG2}
K_G=\frac{\pi}{2\ln(1+\sqrt{2})}. 	
\end{equation}
The conjecture was refuted in 2011 by Braverman, Makarychev and Naor, who showed that $K_R$ is strictly less than Krivine's bound \cite{Braverman2011,Braverman2013}.

Similarly, if the numbers $s_i$ and $t_j$ and matrix ${\bf A}$ are taken as complex, then a similar set of constants $K_C(n)$ may be defined \cite{Krivine1977,Krivine1979,Haagerup1987,Finch2003}, where the upper limit is given by $8/[\pi(x_0+1)]$, $x_0$ being the root of
\begin{equation}
\frac{1}{8}\pi(x+1)=x\int_0^{\pi/2}\frac{\cos^2\theta}{\sqrt{1-x^2\sin^2\theta}}d\theta=\frac{1}{x}\left[E(x)-(1-x^2)K(x)\right],	
\end{equation}
where 
\begin{equation}
K(k)=\int_0^{\pi/2}\frac{d\theta}{\sqrt{1-k^2\sin^2\theta}}
\end{equation}
and
\begin{equation}
E(k)=\int_0^{\pi/2}\sqrt{1-k^2\sin^2\theta}d\theta
\end{equation}
are complete elliptic integrals of the first and second kind respectively. However, Haagerup has suggested that the upper limit would more plausibly be given by \cite{Haagerup1987}
\begin{equation}
\left[\int_0^{\pi/2}\frac{\cos^2\theta}{\sqrt{1+\sin^2\theta}}d\theta\right]^{-1}=\frac{1}{2K(i)-E(i)}.	
\end{equation}
It is worth mentioning that $K_G$ is related to Khintchine's constant \cite{Khintchine1964,Haagerup1982,Jameson1985,Bailey1997}\footnote{If
\begin{equation}
x=[a_0;a_1,...]=a_0+\frac{1}{a_1+\frac{1}{\displaystyle a_2+\frac{1}{a_3+...}}} 	
\end{equation}
is the simple continued fraction of a real number $x$, where the numbers $a_i$ are the partial denominators, Khintchine has shown, using the Gauss-Kuz'min distribution, that for almost all positive irrationals the limiting geometric mean of the positive elements $a_i$ exists and is equal to \cite{Khintchine1964}:
\begin{equation}
K=\lim_{n\rightarrow\infty}\left(a_1a_2...a_n\right)^{1/n}.
\end{equation}
In particular, one has
\begin{equation}
K=\prod_{n=1}^{\infty}\left[1+\frac{1}{n(n+2)}\right]^{\frac{\ln n}{\ln 2}}. 	
\end{equation}
}.
Equation (\ref{KG2}) can be put in the form
\begin{equation}
\ln^2\left(\sqrt{2}+1\right)=\frac{\pi^2}{4K_G^2}
\end{equation}
yielding
\begin{empheq}[box=\fbox]{align}
\sum_{n=1}^{\infty}\left[\sum_{k=n+1}^{\infty}(-1)^{k}\left(\frac{1}{4k-1}-\frac{1}{4k-3}\right)\right]^2=\frac{\pi}{16}\left(1-\frac{\pi}{K_G^2}\right),
\end{empheq}
or
\begin{empheq}[box=\fbox]{align}
K_G=\frac{\pi}{\displaystyle\sqrt{\pi-16\sum_{n=1}^{\infty}\left[\sum_{k=n+1}^{\infty}(-1)^{k}\left(\frac{1}{4k-1}-\frac{1}{4k-3}\right)\right]^2}}.
\end{empheq}

\vspace{1cm}

In 1935, Einstein, Podolsky and Rosen published a famous article suggesting that the quantum mechanical description of reality is incomplete and proposing the existence of hidden variables that cannot be measured, but explain the statistical character of the (experimentally confirmed) predictions of quantum mechanics. The concept of ``locality'' is connected to the assumption that the observations are independent; it constitutes the cornerstone of the Bohm theory involving non-local hidden variables in order to reconcile theory with experiments. In other words, standard quantum mechanics would rely on the statistical description of the underlying hidden variables. In 1964, Bell observed that such an assumption could be tested and proposed an inequality that should be satisfied by those hidden variables. Clauser, Horne, Shimony and Holt \cite{Clauser1969} modified the Bell inequality and suggested that it could be checked experimentally. Many experiments later, there is no doubt that the Bell-Clauser-Horne-Shimony-Holt inequality \cite{Clauser1969} is violated, invalidating thus the ``hidden variables'' theory \cite{Aspect2000a,Aspect2000b}. In 1980 Tsirelson observed that the Grothendieck constant could be interpreted as an upper bound for a generalized Bell inequality, and that such a violation is closely related to the assertion that the Grothendieck constant is greater than 1. He also found a variant of the Clauser-Horne-Shimony-Holt inequality (now called ``Tsirelson's bound'' \cite{Tsirelson1980}).

\section{Conclusion}

We obtained a double-series formula for $\pi$ using the Fourier series expansion of $1/\cos\left(x/4\right)$ and applying the Parseval-Plancherel identity. We pointed out that such a formula can in turn be interpreted as an expression of the Grothendieck-Krivine constant.

\section{Appendix}

It is worth mentioning that it is possible to derive further relations for the Grothendieck constant using the following formulas recently published by Valdebenito \cite{Valdebenito2021}:

\begin{equation}\label{eq1}
\pi^2=4\left[\ln(1+\sqrt{2})\right]^2+8\int_{\ln(1+\sqrt{2})}^{\infty}\arcsinh\left(\csch x\right)dx,
\end{equation}
\begin{equation}
\pi^2+4\left[\ln(1+\sqrt{2})\right]^2+8\int_0^{\ln(1+\sqrt{2})}\arcsinh\left(\csch x\right)dx,
\end{equation}
\begin{equation}
\pi^2=4\left[\ln(1+\sqrt{2})\right]^2+8\int_{\ln(1+\sqrt{2})}^{\infty}\arccosh\left(\coth x\right)dx
\end{equation}
and

\begin{equation}\label{eq4}
\pi^2=4\left[\ln(1+\sqrt{2})\right]^2+8\int_{\ln(1+\sqrt{2})}^{\infty}\arctanh\left(\sech x\right)dx,
\end{equation}
where
\begin{equation}
\csch x=\frac{1}{\sinh x}
\end{equation}
is the usual cosecant function and
\begin{equation}
\sech x=\frac{1}{\cosh x}
\end{equation}
the hyperbolic secant. Using the definition of the Grothendieck-Krivine constant (\ref{KG}), Eqs. (\ref{eq1}) to (\ref{eq4}) give respectively
\begin{equation}
K_G=\frac{1}{\displaystyle\sqrt{1-\frac{8}{\pi^2}\int_{\frac{\pi}{2K_G}}^{\infty}\arcsinh\left(\csch x\right)dx}},
\end{equation}
\begin{equation}
K_G=\frac{1}{\displaystyle\sqrt{\frac{8}{\pi^2}\int_0^{\frac{\pi}{2K_G}}\arcsinh\left(\csch x\right)dx-1}},
\end{equation}
\begin{equation}
K_G=\frac{1}{\displaystyle\sqrt{1-\frac{8}{\pi^2}\int_{\frac{\pi}{2K_G}}^{\infty}\arccosh\left(\coth x\right)dx}},
\end{equation}
and
\begin{equation}
K_G=\frac{1}{\displaystyle\sqrt{1-\frac{8}{\pi^2}\int_{\frac{\pi}{2K_G}}^{\infty}\arctanh\left(\sech x\right)dx}}.
\end{equation}


\begin{thebibliography}{99}

\bibitem{Borwein} J. M. Borwein and P. B. Borwein, {\it Pi and the AGM}, John Wiley and Sons, Inc., New York, Toronto, 1986.

\bibitem{Grothendieck1956} A. Grothendieck, {\it R\'esum\'e de la th\'eorie m\'etrique des produits tensoriels topologiques}, Bol.  Soc. Mat. Sao Paulo {\bf 8}, 1-79 (1956).

\bibitem{Pisier1978} G. Pisier, {\it Grothendieck's theorem for non-commutative C$^{*}$-algebras with an appendix on Grothendieck's constant}, J. Funct. Anal. {\bf 2}, 379-415 (1978). 

\bibitem{Bell1964} J. S. Bell, {\it On the Einstein Podolsky Rosen Paradox}, Physics Physique Fizika {\bf 1}, 195-200 (1964).

\bibitem{Einstein1935} A. Einstein, B. Podolsky and N. Rosen, {\it Can Quantum-Mechanical Description of Physical Reality Be Considered Complete?}, Phys. Rev. {\bf 47}, 777-780 (1935).

\bibitem{Aspect2000a} A. Aspect, {\it Bell's theorem: the naive view of an experimentalist}, Quantum [Un]speakables (Vienna, 2000), 119-153, Springer, Berlin, 2002, (arXiv:quant-ph/0402001).

\bibitem{Aspect2000b} A. Aspect, {\it Testing Bell's inequalities}, Quantum Reflections, 69-78, Cambridge Univ. Press, Cambridge, 2000.

\bibitem{Tsirelson1980} B. S. Tsirelson, {\it Quantum generalizations of Bell's inequality}, Lett. Math. Phys. {\bf 4}, 93-100 (1980).

\bibitem{Tsirelson1993} B. S. Tsirelson, {\it Some results and problems on quantum Bell-type inequalities}, Hadronic J. Suppl. {\bf 8}, 329-345 (1993).

\bibitem{Pisier2012} G. Pisier, {\it Grothendieck theorem, past and present}, Bull. Amer. Math. Soc. {\bf 49}, 237-323 (2012).

\bibitem{Krivine1977} J.-L. Krivine, {\it Sur la constante de Grothendieck}, C. R. A. S. {\bf 284}, 445-446 (1977).

\bibitem{Lelionnais1983} F. Le Lionnais, {\it Les nombres remarquables}, Paris, Hermann, p. 42, 1983.

\bibitem{Braverman2011} M. Braverman, K. Makarychev, Y. Makarychev and A. Naor, {\it The Grothendieck constant is strictly smaller than Krivine's bound}, IEEE 52$^{nd}$ Annual Symposium on Foundations of Computer Science, pp. 453-462 (2011).

\bibitem{Braverman2013} M. Braverman, K. Makarychev, Y. Makarychev and A. Naor, {\it The Grothendieck constant is strictly smaller than Krivine's bound}, Forum of Mathematics, Pi, 1, E4. doi:10.1017/fmp.2013.4 (2013).

\bibitem{Tolstov1962} G. P. Tolstov, {\it Fourier Series}, translated by R. Silverman, Englewood Cliffs, NJ: Prentice-Hall, Inc., 1962.

\bibitem{Kaplan1992} Kaplan, W. {\it Advanced Calculus}, 4$^{th}$ ed., Reading, MA: Addison-Wesley, p. 501, 1992.

\bibitem{Tissier1991} A. Tissier, {\it Math\'ematiques g\'en\'erales - Agr\'egation interne de Math\'ematiques}, Br\'eal, Montreuil, 1991.

\bibitem{Parseval1806} M.-A. Parseval des Ch\^enes, {\it M\'emoire sur les s\'eries et sur l'int\'egration compl\`ete d'une \'equation aux diff\'erences partielles lin\'eaire du second ordre, \`a coefficients constants} in ``M\'emoires pr\'esent\'es \`a l'Institut des Sciences, Lettres et Arts, par divers savants, et lus dans ses assembl\'ees'', Sciences, math\'ematiques et physiques (Savants \'etrangers) {\bf 1}, 638-648 (1806).

\bibitem{Plancherel1910} M. Plancherel, {\it Contribution \`a l'etude de la representation d'une fonction arbitraire par les integrales d\'efinies}, Rendiconti del Circolo Matematico di Palermo {\bf 30}, 298-335 (1910).

\bibitem{mathe} \url{https://mathworld.wolfram.com/GrothendiecksConstant.html}.

\bibitem{Krivine1979} J.-L. Krivine, {\it Constantes de Grothendieck et fonctions de type positif sur les sph\`eres}, Adv. Math. {\bf 31}, 16-30 (1979). 

\bibitem{Haagerup1987} U. Haagerup, {\it A New Upper Bound for the Complex Grothendieck Constant}, Israeli J. Math. {\bf 60}, 199-224 (1987).

\bibitem{Finch2003} S. R. Finch, {\it Grothendieck's Constants}, \S 3.11 in Mathematical Constants. Cambridge, England: Cambridge University Press, pp. 235-237, 2003.

\bibitem{Khintchine1964} A. Khintchine, {\it Continued fractions}, University of Chicago Press, Chicago, 1964.

\bibitem{Haagerup1982} U. Haagerup, {\it The best constants in the Khintchine inequality}, Studia Math. {\bf 70}, 231-283 (1982).

\bibitem{Jameson1985} G. J. O. Jameson, {\it The interpolation proof of Grothendieck's inequality}, Proc. Edinburgh Math. Soc. {\bf 28}, 217-223 (1985).

\bibitem{Bailey1997} D. H. Bailey, J. M. Borwein and R. E. Crandall, {\it On the Khintchine Constant}, Math. Comput. {\bf 66}, 417-431 (1997).

\bibitem{Clauser1969} J. F. Clauser, M. A. Horne, A. Shimony, R. A. Holt, {\it Proposed experiment to test local hidden-variable theories}, Phys. Rev. Lett. {\bf 23}, 880-884 (1969).

\bibitem{Valdebenito2021} E. Valdebenito, {\it Elementary identities and Pi},\\ \url{https://vixra.org/author/edgar_valdebenito}.

\end{thebibliography}
\end{document}